\documentclass[oneside,reqno]{amsart}
\usepackage{amsmath, amsthm, mathrsfs, amssymb, amsfonts, mathtools}
\usepackage{graphicx}
\usepackage{hyperref}
\hypersetup{colorlinks=true, linkcolor=blue, citecolor=magenta, filecolor=magenta, urlcolor=cyan}

\usepackage{doi}

\usepackage[backend=biber, style=trad-alpha, doi=true, isbn=false, hyperref=true, backref=false, giveninits=true, abbreviate=true, sorting=nyt]{biblatex}
\usepackage[dvipsnames]{xcolor}
\newtheorem{remark}{Remark}

\def\diff{\mathrm{d}}
\newcommand{\dd}{\mathrm{d}}
\newcommand{\sgn}{\mathrm{sgn}}
\DeclareMathOperator{\cn}{cn}

\definecolor{verylightgray}{rgb}{0.9,0.9,0.9}

\usepackage{orcidlink}

\usepackage{afterpage}

\numberwithin{equation}{section}
\mathtoolsset{showonlyrefs}

\begin{document}

\title[Structure of time-periodic solutions decoded in Poincar\'{e}-Lindstedt series]{Complex structure of time-periodic solutions decoded in Poincar\'{e}-Lindstedt series: the cubic conformal wave equation on $\mathbb{S}^{3}$}

\author{Filip Ficek\hspace{0.233ex}\orcidlink{0000-0001-5885-7064}}
\address{(F.F.) University of Vienna, Faculty of Mathematics, Oskar-Morgenstern-Platz 1, 1090 Vienna, Austria and University of Vienna, Gravitational Physics, Boltzmanngasse 5, 1090 Vienna, Austria}
\email[]{filip.ficek@univie.ac.at}
\author{Maciej Maliborski\hspace{0.233ex}\orcidlink{0000-0002-8621-9761}}
\address{(M.M.) University of Vienna, Faculty of Mathematics, Oskar-Morgenstern-Platz 1, 1090 Vienna, Austria and University of Vienna, Gravitational Physics, Boltzmanngasse 5, 1090 Vienna, Austria}
\email[]{maciej.maliborski@univie.ac.at}
\thanks{We acknowledge the support of the Austrian Science Fund (FWF) through Project \href{http://doi.org/10.55776/P36455}{P 36455}, Wittgenstein Award \href{http://doi.org/10.55776/Z387}{Z 387}, and the START-Project \href{http://doi.org/10.55776/Y963}{Y 963}. Computations were performed on a supercomputer hosted by the Institute of Theoretical Physics at Jagiellonian University, with funding provided by the Polish National Science Centre grant No. 2017/26/A/ST2/00530.}


\begin{abstract}
This work explores the rich structure of spherically symmetric time-periodic solutions of the cubic conformal wave equation on $\mathbb{S}^{3}$.
We discover that the families of solutions bifurcating from the eigenmodes of the linearised equation form patterns similar to the ones observed for the cubic wave equation.
Alongside the Galerkin approaches, we study them using the new method based on the Pad\'{e} approximants.
To do so, we provide a rigorous perturbative construction of solutions.
Due to the conformal symmetry, the solutions presented in this work serve as examples of large time-periodic solutions of the conformally coupled scalar field on the anti-de Sitter background.
\end{abstract}

\date{\today}

\maketitle

\tableofcontents

\section{Introduction}
\label{sec:Introduction}

\subsection{Model}
\label{sec:Model}

We consider the conformally invariant cubic wave equation
\begin{equation}
    \label{eq:25.03.20_02}
    g^{\mu\nu}\nabla_{\mu}\nabla_{\nu}\phi - \frac{1}{6} R(g)\phi - \lambda\phi^3 = 0
\end{equation}
on the Einstein cylinder $\mathbb{R} \times \mathbb{S}^3$ with the metric
\begin{equation}
    \label{eq:25.03.20_03}
    g = -\diff{t}^2 + \rho^2 \left( \diff{x}^2 + \sin^{2}{x} \,\diff{\Omega}^2 \right)\,,
\end{equation}
where $\rho$ is the radius of the 3-sphere, $\diff\Omega^2$ is the round metric on the unit 2-sphere, $R(g) = 6/\rho^2$ is the Ricci scalar of $g$, and $\lambda$ is a positive constant.

After rescaling $t \to t/\rho$ and $\phi \to \rho \sqrt{\lambda} \phi$, equation \eqref{eq:25.03.20_02} takes the dimensionless form of the cubic Klein-Gordon equation with unit mass
\begin{equation}
    \label{eq:25.03.20_04}
    \phi_{tt} - \Delta_{\mathbb{S}^3} \phi + \phi + \phi^3 = 0\,.
\end{equation}
Assuming spherical symmetry $\phi=\phi(t,x)$ and substituting $\phi(t,x) = u(t, x)/\sin x$ into \eqref{eq:25.03.20_04} we obtain a one-dimensional wave equation with the Dirichlet boundary conditions that are enforced by the regularity of $\phi$ on $\mathbb{S}^3$
\begin{equation}
    \label{eq:25.03.20_05}
    \partial_{t}^{2}u - \partial_{x}^{2}u + \frac{u^{3}}{\sin^{2}{x}} = 0\,,
    \quad
    u(t,0) = 0 = u(t,\pi)\,.
\end{equation}

Our goal is to construct time-periodic solutions of equation \eqref{eq:25.03.20_05}, i.e., solutions satisfying, for some period $T\geq 0$,
\begin{equation}
    u(t+T, x) = u(t, x)\,.
\end{equation}
and study their properties. For the following discussion, it is convenient to rescale time $\tau=\Omega t$, so the period is equal to $2\pi$, $u(\tau,x)=u(\tau+2\pi,x)$, and the equation becomes:
\begin{align}
    \label{eq:conformal2}
    \Omega^2 \partial_\tau^2 u -\partial_x^2 u + \frac{u^3}{\sin^2 x}=0\,.
\end{align}
The conserved energy of \eqref{eq:conformal2} takes the form
\begin{equation}
    \label{eq:25.03.24_05}
    E[u] = \int_{0}^{\pi}\left(\frac{1}{2}\Omega^{2}\left(\partial_{\tau}u\right)^{2}+\frac{1}{2}\left(\partial_{x}u\right)^{2}+\frac{1}{4}\frac{u^{4}}{\sin^{2}{x}}\right)\diff{x}\,.
\end{equation}

It is worth noting that equation \eqref{eq:25.03.20_05} admits a breathing mode (see \cite{Evnin.2020fys})
\begin{equation}
    \label{eq:25.07.23}
    B[u] = \int_{0}^{\pi} \left[\cos{x}\left((\partial_{t}u)^{2} + (\partial_{x}u)^{2} + \frac{1}{2}\frac{u^{4}}{\sin^{2}{x}}\right) - 2i\sin{x}\,\partial_{t}u\,\partial_{x}u\right] \diff{x}\,,
\end{equation}
whose absolute value represents a conserved quantity. However, for the class of solutions considered in this work (solutions symmetric or anti-symmetric across $x=\pi/2$, see the remark below), it vanishes identically and therefore does not provide additional information.

\begin{remark}
An alternative approach to derive the same equation would be to consider a conformally coupled self-interacting scalar field in four-dimensional anti-de Sitter (AdS) spacetime. In this case, the metric is given by
\begin{equation}
    \label{eq:25.03.27_02}
    \tilde{g}=\frac{1}{\cos^{2}{x}}\left(-\diff{t}^{2}+\diff{x}^{2}+\sin^2{x}\, \diff{\Omega}^{2}\right)\,,
\end{equation}
where $x$ ranges from $0$ to $\pi/2$. Since $g=\cos^{2}{\!x}\,\tilde{g}$, applying a field redefinition $\tilde{\phi}(t,x)=\phi(t,x)\cos{x}$ transforms the equation between the Einstein cylinder and AdS space, as follows from equation \eqref{eq:25.03.20_02}.

To specify the evolution, a boundary condition must be imposed at the conformal boundary located at $x = \pi/2$. In this case, we impose the Dirichlet, $\phi(t,\pi/2) = 0$, or Neumann condition, $\partial_{x}\phi(t,\pi/2)=0$, which correspond to the typical reflecting boundary conditions used in the AdS context \cite{Bizoń.2014, FFMM24}. This model is a specific case of the full Einstein cylinder model, which can be recovered by extending the domain to $[0,\pi]$ and enforcing the reflection symmetry through the equator $\phi(t, x) = \mp\phi(t, \pi - x)$, with the sign chosen to match the imposed boundary condition.

We highlight this connection because the study of nonlinear dynamics of small perturbations in AdS has received considerable interest in recent years \cite{Evnin.2021}, and the objectives of this work align with this area of research.
\end{remark}

\begin{remark}
    In this work, we focus exclusively on the defocusing nonlinearity. While the techniques used here can be easily adapted to the focusing case, the results exhibit significant differences. We postpone their detailed analysis for a future study.
\end{remark}

\subsection{Related studies}
\label{sec:RelatedStudies}

Time periodic solutions of equation \eqref{eq:25.03.20_04} with quadratic, cubic, and quintic nonlinearity have recently been an object of intense study \cite{Bizon.2017,  Chatzikaleas.2020oh, Berti.2024,Silimbani.2024, Chatzikaleas.2024}.
In particular, for the cubic nonlinearity, the authors of \cite{Chatzikaleas.2024} used the result of Bambusi–Paleari \cite{Bambusi.2001} to prove the existence of small amplitude spherically symmetric solutions bifurcating from a single eigenmode of the linearized problem ($\sin{nx}$, $n\in\mathbb{N}_{+}$). This approach relies on the stability of stationary solutions of the resonant approximation that were earlier studied in \cite{Bizoń.2017ndf,Bizoń.2019,Bizoń.2020qd}.
In addition, \cite{Chatzikaleas.2024} proves the existence of a finite number of families of aspherical solutions arising from a special ansatz for Hopf fibration of $\mathbb{S}^3$ (the so-called Hopf-plane waves) restricted to particular momenta and bifurcating from a single eigenmode; for details see \cite{Evnin.2021,Chatzikaleas.2024} and references therein.
To avoid the small-divisor problem \cite{arnold1988dynamical} the proof assumes the set of frequencies satisfies a certain Diophantine condition, introduced in \cite{Bambusi.2001}, which implies it is uncountable, accumulates at one, and is of measure zero.
Using variational methods of mountain pass type, the work \cite{Berti.2024} generalises \cite{Chatzikaleas.2024} to quadratic and quintic nonlinearity in spherical symmetry and proves the existence and multiplicity of time-periodic solutions whose frequency again belongs to a Cantor set accumulating at one with measure zero. Contrary to the cubic case, these solutions do not bifurcate from a single mode, but rather from a non-trivial combination, c.f. \cite{Ficek.2024A} and the discussion below.
In addition, for cubic nonlinearity \cite{Berti.2024} extends the result of \cite{Chatzikaleas.2024} to any momenta of the Hopf-plane waves.
The following work \cite{Silimbani.2024} proved the existence and multiplicity of positive measure Cantor families of small amplitude solutions with cubic nonlinearity. To enlarge the set of frequencies for which solutions do exist, this work used the Lyapunov-Schmidt decomposition and a Nash-Moser iteration.

The cited works represent the current state of research on time-periodic solutions to resonant Hamiltonian PDEs in dimensions greater than one.
There is an extensive body of literature on semilinear wave and Klein-Gordon equations in one dimension; see e.g.\ \cite{berti2007nonlinear} and references therein. Among these types of problems, significant interest has been attracted by the one-dimensional cubic wave equation with Dirichlet boundary conditions
\begin{equation}
    \label{eq:25.03.20_01}
    \partial_{t}^{2}u - \partial_{x}^{2}u + u^{3} = 0\,,
    \quad
    u(t,0) = 0 = u(t,\pi)\,.
\end{equation}
Existence of certain time-periodic solutions in this setup has been proven in \cite{Bambusi.2001, Berti.2003, Berti.2004}.
However, those results rely on the local bifurcation theory and necessarily focus on small solutions.
They also have to deal with the small divisor problem, leading to nowhere dense frequency sets. In \cite{Ficek.2024A, Ficek.2024B} we have described an intricate pattern that is formed by large time-periodic solutions to equation \eqref{eq:25.03.20_01}.
We have also hinted that similar structures can be observed in numerous other nonlinear dispersive equations.

The goal of this work is to investigate such structures in a specific model with appealing geometric properties, serving as a toy model related to General Relativity with a negative cosmological constant \cite{Bizoń.2017ndf}.
Consequently, this work complements previous existence proofs of time-periodic solutions for completely resonant Hamiltonian PDEs on manifolds with dimension greater than one.

\subsection{Features of the model}
\label{sec:FeaturesOfTheModel}

There are some important differences between the conformal equation discussed in this work and the cubic equation \eqref{eq:25.03.20_01} studied in \cite{Ficek.2024A, Ficek.2024B}. Most importantly, here the spatial modes get mixed by the nonlinearity according to the formula
\begin{align}
    \label{eq:decomposition_main}
    \frac{\sin jx\, \sin k x\, \sin l x }{\sin^2 x}=\sum_{m=1}^{j+k+l-2} S_{jklm}\, \sin mx\,,
\end{align}
where $S_{jklm}$ are the interaction coefficients. Both Eq.~\eqref{eq:decomposition_main} and explicit formulas for $S_{jklm}$ are derived in the Appendix~\ref{sec:InteractionCoefficients}. 
This formula strongly differs from the one characteristic for the cubic nonlinearity, where interactions are much sparser and the non-zero interaction coefficients are simpler. This difference makes the approach based on the reducible systems (see \cite{Ficek.2024B}) more complicated for the conformal equation. However, there is also an important advantage. The highest mode in decomposition \eqref{eq:decomposition_main} is $j+k+l-2$, while for the simple cubic nonlinearity \eqref{eq:25.03.20_01}, an analogous formula would give $j+k+l$.
This difference could be attributed to the high symmetry of the AdS spacetime, as explained in \cite{Evnin.2016}, which due to the conformal relation mentioned in the Remark~1 is also relevant in our case.
As a result, we are able to give a perturbative construction of a family of solutions bifurcating from a single mode. For Eq.~\eqref{eq:25.03.20_01}, such a feat is impossible, since already at the linear level one needs a highly non-trivial infinite combination of modes in order to remove all the resonances in the next order, see \cite{Ficek.2024A}.

Another consequence of the presence of expression $\sin^2 x$ in the nonlinearity is the lack of the scaling symmetry. In the case of the cubic wave equation \eqref{eq:25.03.20_01} from any solution $u(t,x)$ with period $T$ we can obtain new solutions $n\,u(nt,nx)$ with periods $T/n$, where $n\in\mathbb{N}_+$. For this reason \cite{Ficek.2024A, Ficek.2024B} focus on a single solution family bifurcating at $\Omega=1$ as other families can be obtained via simple rescalings. This is not the case for the conformal equation \eqref{eq:25.03.20_05} since the nonlinear term breaks the mentioned symmetry. Hence, we have to treat solution families bifurcating from different eigenmodes of the linearised equation separately. In this regard, the solution family bifurcating from the eigenmode $\sin x$ at $\Omega=1$ is special since it can be described by an explicit formula. The next subsection is devoted to this subject. The remaining part of the paper discusses solutions bifurcating from $\sin N x$ at $\Omega=N$, where $N\geq2$.

\subsection{Exact solution}
\label{sec:ExactSolution}

Before considering the more general case, we derive an explicit expression for the solution bifurcating from the lowest eigenmode $\sin{x}$. For the single mode ansatz $u(t, x) = \phi(t) \sin x$, equation \eqref{eq:25.03.20_05} reduces to the undamped, unforced Duffing equation 
\begin{equation}
    \label{eq:25.03.24_01}
    \ddot{\phi} + \phi + \phi^3 = 0\,,
\end{equation}
which is just \eqref{eq:25.03.20_04} for solutions that are homogeneous in space. For initial data $(\phi,\dot{\phi})|_{t=0}=(\varepsilon,0)$ with any $\varepsilon$, this equation is solved by
the Jacobi elliptic function \cite[\href{https://dlmf.nist.gov/22.2}{Sec.~22.2}]{DLMF}
\begin{equation}
    \label{eq:25.03.24_02}
    \phi(t; \varepsilon) = \varepsilon \cn\left(\sqrt{1+\varepsilon^2}\,t, \kappa\right), \quad \kappa(\varepsilon) = \frac{\varepsilon}{\sqrt{2(1+\varepsilon^2)}}\,.
\end{equation}
Since the function $\cn(t,\kappa)$ has period $4K(\kappa)$, where $K(\kappa)$ is the complete elliptic integral of the first kind \cite[\href{https://dlmf.nist.gov/22.2}{Sec.~22.2}]{DLMF}
\begin{equation}
    \label{eq:25.03.24_03}
    K(\kappa) := \int_0^{\pi/2} \frac{d\theta}{\sqrt{1 - \kappa^2 \sin^2\theta}} = \frac{\pi}{2} {}_2F_1\left(\frac{1}{2}, \frac{1}{2}, 1; \kappa^2\right)\,,
\end{equation}
the solution $\phi(t;\varepsilon)$ has frequency
\begin{equation}
    \label{eq:25.03.24_04}
    \Omega(\varepsilon) = \frac{\pi}{2}\frac{\sqrt{1+\varepsilon^{2}}}{K(\kappa)}
    \,.
\end{equation}
The frequency $\Omega(\varepsilon)$ is a monotone increasing function of $\varepsilon$ and behaves as $\Omega(\varepsilon)=1+3\varepsilon^{2}/8 + \mathcal{O}(\varepsilon^{4})$ for small $\varepsilon$ and as $\Omega(\varepsilon)=\sqrt{2\pi}\varepsilon\Gamma(3/4)/\Gamma(1/4) + \mathcal{O}(\varepsilon^{-1})$ for $\varepsilon\to\infty$.
The energy associated with the solution given by Eq.~\eqref{eq:25.03.24_02} is $E(\varepsilon) = \pi\varepsilon^{2}(2+\varepsilon^{2})/8$.
A perturbative construction of the single-mode configuration was presented in \cite{Chatzikaleas.2020oh}.
Clearly, equation \eqref{eq:25.03.20_01} has no single-mode solutions.

\subsection{Main results}
\label{sec:MainResults}

This work provides a first report on and examines the intricate structure of time-periodic solutions to the cubic conformal wave equation on $\mathbb{S}^{3}$ (being an example of a Hamiltonian PDE on a manifold of dimension larger than one). Using a Galerkin scheme, we analyse its lowest-dimensional truncations to discover the bifurcation structure of solutions. By studying how this structure changes as the truncation increases we hypothesise about the solutions to the PDE.

Our results suggest existence of a rich families of solutions to \eqref{eq:conformal2}.
From every linear eigenmode $\sin{Nx}$, $N\in\mathbb{N}_{+}$, there bifurcates a solution curve that consists of a main part that we call a trunk.
(The part of the trunk closest to the bifurcation point contains solutions previously described in the literature.)
Furthermore, for $N>1$ we identify the emergence of branches from the respective trunks, a phenomenon similar to the one observed and described in \cite{Ficek.2024A, Ficek.2024B}. Those branches densely populate each trunk and their existence correlates with Cantor sets appearing in rigorous existence proofs \cite{Berti.2024,Silimbani.2024, Chatzikaleas.2020oh, Chatzikaleas.2024}.
For the lowest eigenmode $N=1$, no branches are observed, and this special solution was given explicitly above.

Additionally, for every $N>1$ we construct a countable family of solutions bifurcating from the linear eigenmodes through a rigorous perturbative expansion. We prove that this expansion extends to an arbitrary order, providing explicit formulas for the Fourier coefficients. This generalises \cite{Chatzikaleas.2020oh} to solutions bifurcating from an arbitrary mode.
To enhance the accuracy of our analysis, we employ a high-precision pseudo-spectral numerical scheme for evaluating projections of the nonlinear terms, enabling us to compute perturbative expansions of time-periodic solutions to very high orders.

A key contribution of this study is the use of Pad\'{e} approximants to the Poincar\'{e}-Lindstedt series. We demonstrate that Pad\'{e} approximants not only successfully recreate the shape of the main trunk (see Fig.~\ref{fig:pade1}) but also encode information about the locations of the branches, cf. Fig.~\ref{fig:spectrum}. Results of this approach constitute an independent confirmation of the numerical findings about the structure of time-periodic solutions to \eqref{eq:conformal2}.

Finally, exploiting the conformal symmetry of the equation, we establish that the solutions constructed in this work serve as explicit examples of large time-periodic solutions of the conformally coupled scalar field on the AdS background with respective boundary conditions (Dirichlet for $N$ being even and Neumann for $N$ odd, respectively). This connection highlights the broader significance of our results within the study of nonlinear wave dynamics in curved spacetime geometries.

\subsection{Structure of the paper}
\label{sec:StructureOfThePaper}

The remainder of this paper is organized as follows.
In Sec.~\ref{sec:PerturbativeExpansion} we develop a perturbative expansion for time-periodic solutions, which serves as the basis for the subsequent analysis of Pad\'{e} approximants.
In Sec.~\ref{sec:NumericalApproaches} we present modifications to the Galerkin method proposed in \cite{Ficek.2024A}, adapting the numerical scheme to Eq.~\eqref{eq:25.03.20_05}.
We also comment on how the numerical method can be used to significantly accelerate the construction of high-order perturbative series.
Sec.~\ref{sec:StructureOfSolutions} is devoted to a detailed examination of the structure of time-periodic solutions. We analyse both numerical and perturbative results, showing how certain key features of solutions can be effectively captured using reducible systems.
Finally, we calculate Pad\'{e} approximants and compare their predictions with numerical data, demonstrating that the perturbative approach encodes essential aspects of the bifurcation structure, such as the emergence and locations of branches.

\section{Perturbative expansion}
\label{sec:PerturbativeExpansion}

As a first step in the construction of the Poincar\'{e}-Lindstedt series, we redefine $u\to\sqrt{\varepsilon}u$ so that Eq.~\eqref{eq:conformal2} becomes
\begin{align}
    \label{eq:conformal_epsilon}
    \Omega^2 \partial_\tau^2 u -\partial_x^2 u + \varepsilon\frac{u^3}{\sin^2 x}=0\,.
\end{align}
Then we expand both the square of frequency $\Omega^2$ and solution $u$ into formal series in $\varepsilon$
\begin{align}
    \label{eq:1_order}
    \Omega^2=N^2+\sum_{n=1}^\infty \varepsilon^n \omega^{(n)}\,,
    \quad
    u(\tau, x)= \sum_{n=0}^\infty \varepsilon^n u^{(n)} (\tau,x)\,,
\end{align}
with $N\geq 1$ and where $u^{(n)}$ are functions periodic in $\tau$ and satisfying Dirichlet boundary conditions in $x$. The small parameter $\varepsilon$ is fixed by
\begin{equation}
    \label{eq:epsilon}
    u^{(0)}(\tau,x)=\cos\tau \, \sin Nx\,,
    \quad
    \int_0^\pi u^{(n)}(0,x) \, \sin Nx \, \dd x = 0 \,\ \mbox{ for } n\geq 1\,.
\end{equation}
Plugging expressions \eqref{eq:1_order} into Eq.~\eqref{eq:conformal2} and expanding around zero, we obtain a  hierarchy of equations enumerated by the increasing order of $\varepsilon$. We will show that it is possible to construct a formal solution in such a way that at every perturbative order the secular terms vanish and the resulting solution is time periodic.\footnote{At this point we treat $\Omega^2$ and $u$ as formal series in $\varepsilon$ and to not consider their convergence.}

The lowest (zeroth) order equation is 
\begin{align*}
    N^2\partial_\tau^2 u^{(0)} -\partial_x^2 u^{(0)}=0\,.
\end{align*}
It is linear and is simply solved by $u^{(0)}$ defined above. Let us point out that such an equation is solved by any function being a sum of terms of type $\cos n\tau\, \sin nNx$, where $n\in\mathbb{N}_+$. We will call such terms resonant. For future reference, let us recall that if $\cos j \tau \, \sin k x$ ($j,k\in\mathbb{N}_+$) is not a resonant term, the equation
\begin{align*}
    N^2 \partial_\tau^2 \psi -\partial_x^2 \psi = f_{jk}  \cos j \tau \, \sin k x\,.
\end{align*}
has a solution in the form of
\begin{align}\label{eq:nonresonant_solution}
    \psi(\tau,x)=\sum_{m=1}^{\infty} a_{m} \cos m \tau\,\sin m N x + \frac{f_{jk}}{k^2-j^2N^2} \cos j \tau\,  \sin k x\,,
\end{align}
where the homogeneous part $a_m$ can be freely chosen. 

In the next order, the equation is
\begin{align}
    \label{eq:1st_order}
    N^2\partial_\tau^2 u^{(1)} -\partial_x^2 u^{(1)}=-\omega^{(1)} \partial_\tau^2 u^{(0)} -\frac{\left(u^{(0)}\right)^3}{\sin^2 x}.
\end{align}
Such an equation has a time-periodic solution satisfying the Dirichlet boundary conditions only if its right-hand side contains no resonant terms. With $u^{(0)}$ fixed as above and using the notation of Appendix~\ref{sec:InteractionCoefficients} for the interaction coefficients $S_{ijkl}$ we can write
\begin{align*}
    -\omega^{(1)} \partial_\tau^2 u^{(0)} -\frac{\left(u^{(0)}\right)^3}{\sin^2 x}=\omega^{(1)} \cos\tau\,\sin Nx-\frac{1}{4}\left(3\cos \tau+\cos 3\tau\right)\sum_{k=1}^{3N-2}S_{NNNk}\sin k x\,.
\end{align*} 
Let us point out that due to the symmetries of $S_{NNNk}$, the sum contains only terms with $k$ of the same parity as $N$. Clearly, the only resonant term that can be present on the right-hand side is $\cos\tau \, \sin Nx$. Let us introduce $\mathbf{Rsn}[f]$ as a notation for the \textit{resonant part of the expression} $f$. Then, the earlier discussion lets us write
\begin{align*}
    \mathbf{Rsn}\left[-\omega^{(1)} \partial_\tau^2 u^{(0)} -\frac{\left(u^{(0)}\right)^3}{\sin^2 x}\right]&=\omega^{(1)} \cos\tau\,\sin Nx-\frac{3}{4}S_{NNNN}\cos \tau\,\sin N x
    \\
    &=\left(\omega^{(1)}-\frac{3N}{4}\right)\cos \tau\,\sin N x\,,
\end{align*} 
where we have used the explicit form of $S_{NNNN}=N$ from the Appendix \ref{sec:InteractionCoefficients}. Thus, by fixing $\omega^{(1)}=3N/4$ we get rid of the resonant terms in the right-hand side of Eq.~\eqref{eq:1st_order}. As a result, we may write the solution to Eq.~\eqref{eq:1st_order} using \eqref{eq:nonresonant_solution}:
\begin{align*}
    u^{(1)}(\tau,x)=a^{(1)}_1 \cos \tau\,\sin Nx + a^{(1)}_3 \cos 3\tau\,\sin 3Nx+\sum_{j=0}^1 \sideset{}{'}\sum_{k=1}^{3N-2}b^{(1)}_{2j+1,k}\cos(2j+1)\tau\,\sin kx\,,
\end{align*}
where the coefficients are given by
\begin{align*}
    b^{(1)}_{1,k}=-\frac{3S_{NNNk}}{4(k^2-N^2)}\,,\qquad b^{(1)}_{3,k}=-\frac{S_{NNNk}}{4(k^2-9N^2)}\,.
\end{align*}
Here and in the following $\sum'$ denotes a sum that omits resonance terms. Let us point out that we have decided to keep two resonant terms with coefficients $a^{(1)}_1$ and $a^{(1)}_3$. The first one is necessary to make sure that \eqref{eq:epsilon} holds at the level of $\varepsilon$. In order to do so we fix
\begin{align*}
    a^{(1)}_1=-b^{(1)}_{3,N}=-\frac{1}{32N}\,.
\end{align*}
The value of $a^{(1)}_3$ will be chosen in the next order in a way that removes the resonant term $\cos3\tau\,\sin3 Nx$.

In the second order, we have
\begin{align}
    \label{eq:2nd_order}
    N^2\partial_\tau^2 u^{(2)} -\partial_x^2 u^{(2)}=-\omega^{(2)} \partial_\tau^2 u^{(0)} -\omega^{(1)} \partial_\tau^2 u^{(1)} -\frac{3\left(u^{(0)}\right)^2 u^{(1)}}{\sin^2 x}.
\end{align}
Let us denote a \textit{non-resonant part of the expression} $f$ as $\mathbf{NRsn}[f]$ so we have a decomposition $f=\mathbf{Rsn}[f]+\mathbf{NRsn}[f]$. While $\mathbf{NRsn}[u^{(1)}]$ is fully determined, there is still some freedom in $\mathbf{Rsn}[u^{(1)}]$ as $a^{(1)}_3$ is not yet fixed. Let us see how it interacts with the nonlinearity by considering
\begin{align*}
    \frac{3\left(u^{(0)}\right)^2\, \mathbf{Rsn}[u^{(1)}]}{\sin^2 x}=&-\frac{3}{32N}\cos^3\tau\frac{\sin^3Nx}{\sin^2 x}+3a^{(1)}_3 \cos^2\tau\cos 3\tau \frac{\sin^2Nx\sin3Nx}{\sin^2 x}
    \\
    =&-\frac{3}{128N}(3\cos\tau+\cos3\tau)\sum_{k=1}^{3N-2}S_{NNNk}\sin k x\\
    &+\frac{3}{4}a^{(1)}_3(\cos\tau+2\cos3\tau+\cos5\tau)\sum_{k=1}^{5N-2}S_{NN,3N,k}\sin k x\,.
\end{align*}
The resonant part of this expression can be easily extracted as
\begin{align*}
    \mathbf{Rsn}\left[\frac{3\left(u^{(0)}\right)^2\, \mathbf{Rsn}[u^{(1)}]}{\sin^2 x}\right]=-\frac{9}{128}\cos\tau\sin Nx+\frac{3N}{2}a^{(1)}_3\cos3\tau\sin 3Nx\,.
\end{align*}
In addition to formulas for the interaction coefficients, we have also used here the fact that $S_{NN,3N,N}=0$, as discussed in Appendix \ref{sec:InteractionCoefficients}. Now we can use these calculations to write the right-hand side of Eq.~\eqref{eq:2nd_order} as
\begin{multline*}
    \mathbf{Rsn}\left[-\omega^{(2)} \partial_\tau^2 u^{(0)} -\omega^{(1)} \partial_\tau^2 u^{(1)}-\frac{3\left(u^{(0)}\right)^2 u^{(1)}}{\sin^2 x}\right] =
    \\
    =\left(\omega^{(2)}-\frac{3}{128}\right)\cos \tau\, \sin N x + \frac{27N}{4} a_{3}^{(1)}\cos 3 \tau\, \sin 3N x- \mathbf{Rsn}\left[\frac{3\left(u^{(0)}\right)^2 u^{(1)}}{\sin^2 x}\right]
    \\
    = \left(\omega^{(2)}+\frac{3}{64}\right)\cos \tau\, \sin N x + \frac{21N}{4} a_{3}^{(1)}\cos 3 \tau\, \sin 3N x- \mathbf{Rsn}\left[\frac{3\left(u^{(0)}\right)^2\, \mathbf{NRsn}[u^{(1)}]}{\sin^2 x}\right]\, ,
\end{multline*}
An important observation is that since $\mathbf{NRsn}[u^{(1)}]$ contains in its spatial part no modes higher than $\sin(3N-2)x$, the only resonant terms that may be present in the right-hand side of Eq.~\eqref{eq:2nd_order} are $\cos \tau\, \sin N x$ and $\cos 3 \tau\, \sin 3N x$. They can be removed by a suitable choice of $\omega^{(2)}$ and $a_3^{(1)}$, respectively. Thus, we fix
\begin{align*}
    \omega^{(2)}= -\frac{3}{64}+\left[\frac{3\left(u^{(0)}\right)^2\, \mathbf{NRsn}[u^{(1)}]}{\sin^2 x}\right]_{1,N}\,,
    \qquad
    a^{(1)}_3=\frac{4}{21N}\left[\frac{3\left(u^{(0)}\right)^2\, \mathbf{NRsn}[u^{(1)}]}{\sin^2 x}\right]_{3,3N}\,,
\end{align*}
where $[f]_{j,k}$ denotes a \textit{coefficient of the $\cos j\tau\,\sin kx$ term in the expression $f$}. Since now the right-hand side of Eq.~\eqref{eq:2nd_order} is fully determined and contains only non-resonant terms, we can write the solution using \eqref{eq:nonresonant_solution} as
\begin{align*}
    u^{(2)}(\tau,x)=\sum_{m=0}^2 a^{(2)}_{2m+1} \cos(2m+1) \tau\,\sin (2m+1)Nx +\sum_{j=0}^2 \sideset{}{'}\sum_{k=1}^{5N-2}b^{(2)}_{2j+1,k}\cos(2j+1)\tau\,\sin kx\,,
\end{align*}
with
\begin{align*}
    b^{(2)}_{2j+1,k}=\frac{1}{k^2-(2j+1)^2N^2}\left[-\omega^{(2)} \partial_\tau^2 u^{(0)} -\omega^{(1)} \partial_\tau^2 u^{(1)} -\frac{3\left(u^{(0)}\right)^2 u^{(1)}}{\sin^2 x}\right]_{2j+1,k}\,.
\end{align*}
We have also used here the fact that $(u^{(0)})^2\,u^{(1)}/\sin^2 x$ has no modes higher than $\cos 5 \tau$ in the temporal part and $\sin(5N-2)x$ in the spatial one. As previously, we fix
\begin{align*}
    a^{(2)}_1=-b^{(2)}_{3,N}-b^{(2)}_{5,N}
\end{align*}
to satisfy Eq.~\eqref{eq:epsilon}, while $a^{(2)}_3$ and $a^{(2)}_5$ will be chosen in the next order so, together with $\omega^{(3)}$, they remove all resonant terms.

This procedure can be continued up to an arbitrary order of $\varepsilon$ giving us a series of functions  
\begin{multline}
    \label{eq:un}
    u^{(n)}(\tau,x)=\sum_{m=0}^{n}a^{(n)}_{2m+1} \cos (2m+1)\tau\,\sin (2m+1)Nx
    \\
    +\sum_{j=0}^{n} \sideset{}{'}\sum_{k=1}^{(2n+1)N-2}b^{(n)}_{2j+1,k}\cos(2j+1)\tau\,\sin kx\,.
\end{multline}
To justify this claim it is sufficient to show that at each order we are able to choose $\omega^{(n)}$ and $a^{(n-1)}_k$ in such a way that all resonances in Eq.~\eqref{eq:conformal_epsilon} at this order vanish. Then, $b_{(2k+1),j}^{(n)}$ can be calculated using Eq.~\eqref{eq:nonresonant_solution} and $a^{(n)}_1$ is fixed by condition \eqref{eq:epsilon} to be
\begin{align*}
    a^{(n)}_1 = -\sum_{k=1}^n b_{(2k+1),N}^{(n)}\, .
\end{align*}
To see this, let us write \eqref{eq:conformal_epsilon} at order $\varepsilon^n$ as
\begin{align}
    \label{eq:epsn}
    N^2 \partial_\tau^2 u^{(n)} -\partial_x^2 u^{(n)}=-\omega^{(n)} \partial_\tau^2 u^{(0)}-\omega^{(1)} \partial_\tau^2 u^{(n-1)}-\frac{3\left(u^{(0)}\right)^2 u^{(n-1)}}{\sin^2 x}-f^{(n)}\,,
\end{align}
where
\begin{align}
    \label{eq:fn}
    f^{(n)}=\frac{1}{\sin^2 x}\sum_{j=0}^{n-1} \sum_{k=0}^{n-j-1}u^{(j)} u^{(k)} u^{(n-j-k-1)}-3\frac{\left(u^{(0)}\right)^2 u^{(n-1)}}{\sin^2 x}+\sum_{k=1}^{n-2}\omega^{(n-k)}\partial_\tau^2 u^{(k)}\,.
\end{align}
One can easily notice that $f^{(n)}$ is fully established, since it does not contain $\omega^{(n)}$, nor $a^{(n-1)}_{2m+1}$, and it contains terms with temporal modes up to $\cos (2n+1)\tau$ and spatial modes up to $\sin ((2n+1)N-2)x$. 
In particular, there are no resonant terms in $f^{(n)}$ higher than $\cos (2n-1)\tau\,\sin(2n-1)Nx$. Now we can treat the remaining part of the right-hand side of Eq.~\eqref{eq:epsn} as before and decompose its resonant part as
\begin{multline*}
    \mathbf{Rsn}\left[-\omega^{(n)} \partial_\tau^2 u^{(0)}-\omega^{(1)} \partial_\tau^2 u^{(n-1)}-\frac{3\left(u^{(0)}\right)^2 u^{(n-1)}}{\sin^2 x}\right]=
    \\
    \omega^{(n)} \cos \tau\, \sin N x + \frac{3N}{4}\sum_{m=0}^{n-1} a_{2m+1}^{(n-1)}(2m+1)^2\cos(2m+1)\tau\, \sin (2m+1)Nx
    \\
    - \mathbf{Rsn}\left[\frac{3\left(u^{(0)}\right)^2\, \mathbf{Rsn}[u^{(n-1)}]}{\sin^2 x}\right]-\mathbf{Rsn}\left[\frac{3\left(u^{(0)}\right)^2 \,\mathbf{NRsn}[u^{(n-1)}]}{\sin^2 x}\right]\, .
\end{multline*}
Since $\mathbf{NRsn}[u^{(n-1)}]$ is already determined, we focus only on the part containing $\mathbf{Rsn}[u^{(n-1)}]$. It can be rewritten as
\begin{multline*}
    \mathbf{Rsn}\left[\frac{3\left(u^{(0)}\right)^2 \,\mathbf{Rsn}[u^{(n-1)}]}{\sin^2 x}\right]
    \\
    = \sum_{m=0}^{n-1}3a_{2m+1}^{(n-1)} \mathbf{Rsn}\left[\cos^2\tau\,\cos(2m+1)\tau\frac{\sin^2 Nx\, \sin(2m+1)Nx}{\sin^2 x}\right]
    \\
    = \sum_{m=0}^{n-1}\frac{3}{4}a_{2m+1}^{(n-1)} \mathbf{Rsn}\left[\left(\cos(2m+3)\tau+2\cos(2m+1)\tau+\cos(2m-1)\tau\right)\times\right.
    \\
    \left.
    \times\sum_{l=1}^{(2m+3)N-2}S_{NN,(2m+1)N,l}\sin lx\right]\,.
\end{multline*}
Let us point out that for any fixed $m\geq1$ the upper limit of summation results in no resonant term $\cos(2m+3)\tau\,\sin(2m+3)N x$. As $S_{NN,(2m+1)N,(2m-1)N}=0$, there is also no resonant term $\cos(2m-1)\tau\,\sin(2m-1)N x$. Thus, we have
\begin{multline*}
    \mathbf{Rsn}\left[\frac{3\left(u^{(0)}\right)^2 \,\mathbf{Rsn}[u^{(n-1)}]}{\sin^2 x}\right] =
    \\
    = \frac{9N}{2}a_{1}^{(n-1)}\cos\tau\,\sin Nx+\sum_{m=1}^{n-1}\frac{3N}{2}a_{2m+1}^{(n-1)}\cos(2m+1)\tau\,\sin (2m+1)Nx\,.
\end{multline*}
Finally, we write the resonant part of the right-hand side of Eq.~\eqref{eq:epsn} as 
\begin{multline*}
    \mathbf{Rsn}\left[-\omega^{(n)} \partial_\tau^2 u^{(0)}-\omega^{(1)} \partial_\tau^2 u^{(n-1)}-\frac{3\left(u^{(0)}\right)^2 u^{(n-1)}}{\sin^2 x}-f^{(n)}\right] =
    \\
   = \omega^{(n)}\cos\tau\,\sin Nx+ \sum_{m=1}^{n-1}\frac{3N}{4}a_{2m+1}^{(n-1)}\left[(2m+1)^2-2\right]\cos(2m+1)\tau\,\sin (2m+1)Nx\\-\frac{9N}{2}a_{1}^{(n-1)}\cos\tau\,\sin Nx-\textbf{Res}\left[\frac{3\left(u^{(0)}\right)^2 \,\mathbf{NRsn}[u^{(n-1)}]}{\sin^2 x}+f^{(n)}\right]\,.
\end{multline*}
Then it is clear that any resonant term that may result from the last two parts can be removed by a suitable choice of $\omega^{(n)}$, for $\cos\tau\,\sin Nx$ term, and $a^{(n-1)}_{2m+1}$ with $m\geq1$, for $\cos(2m+1)\tau\,\sin (2m+1)Nx$ terms. Afterwards, we can calculate $b^{(n)}_{(2j+1),k}$ and $a^{(n)}_{1}$ and move on to the next order.

In the end, we get the following formulas for coefficients in \eqref{eq:un} with $n\geq 2$:
\begin{subequations}
    \label{eq:25.04.07_01}
\begin{align}
    \label{eq:25.04.07_01a}
    \omega^{(n)}=&\frac{3N}{2}a_{1}^{(n-1)}+\left[3\frac{\left(u^{(0)}\right)^2\, \mathbf{NRsn}[u^{(n-1)}]}{\sin^2 x}+f^{(n)}\right]_{1,N}\,,\\
    \label{eq:25.04.07_01b}
    a_{(2m+1)}^{(n-1)}=&\frac{4}{3N[(2m+1)^2-2]}\left[3\frac{\left(u^{(0)}\right)^2\, \mathbf{NRsn}[u^{(n-1)}]}{\sin^2 x}+f^{(n)}\right]_{(2m+1),(2m+1)N}\,,\\
    &\mbox{for $m\in \{1,2,...,n-1\}$},\\
    \label{eq:25.04.07_01c}
    b_{(2j+1),k}^{(n)}=&\frac{1}{k^2-(2j+1)^2 N^2}\left[-\omega^{(1)} \partial_\tau^2 u^{(n-1)}-3\frac{\left(u^{(0)}\right)^2 u^{(n-1)}}{\sin^2 x}-f^{(n)} \right]_{(2j+1),k}\,,\\
    &\mbox{for $j\in \{1,2,...,n\}$, $k\in\{...,(2n+1)N-4,(2n+1)N-2\}$,}\\
    \label{eq:25.04.07_01d}
    a_{1}^{(n)}=&-\sum_{j=0}^n b_{(2j+1),N}^{(n)}\,,
\end{align}
\end{subequations}
where sums in $j$ are over every second number and start with $j=1$ if $N$ is odd, or with $j=2$ if $N$ is even. 

\section{Numerical approaches}
\label{sec:NumericalApproaches}

\subsection{Galerkin-based scheme}
\label{sec:GalerkinBasedScheme}

To explore the structure of time-periodic solutions beyond the perturbative regime, we use the numerical scheme, based on the Galerkin approach, initially designed for the cubic wave equation \cite{Ficek.2024A}. Here we review this method briefly, as the adaptation to the current problem is rather minor. For improved performance, it requires considering even and odd $N$ cases separately. In the following, we assume $N$ is even. The other case requires a different choice of basis functions and compatible grid points and weights.

We approximate solutions within a finite-dimensional subspace of a Hilbert space
\begin{equation}
	\label{eq:24.05.09_04}
	\textrm{span}\left\{\left.\cos{(2j+1)\tau}\, \sin{2(k+1)x}\,\right|\ j=0,\ldots,M_{\tau}-1\,,\ k=0,\ldots,M_{x}-1\right\}
	\,,
\end{equation}
meaning we write $u(\tau,x)$ as a finite Fourier series
\begin{equation}
	\label{eq:24.05.13_02}
	u_{M_{\tau},M_{x}}(\tau,x) = \sum_{j=0}^{M_{\tau}-1}\sum_{k=0}^{M_{x}-1}\hat{u}_{jk}\cos(2j+1)\tau\, \sin2(k+1)x
	\,.
\end{equation}
In the domain $(\tau,x)\in[0,2\pi]\times[0,\pi]$, we take compatible collocation points
\begin{equation}
	\label{eq:24.05.13_03}
	\tau_{j} = \frac{\pi(j+1/2)}{2M_{\tau}+1}\,,
	\quad
	j=0,\ldots,M_{\tau}-1\,,
	\quad
	x_{k} = \frac{\pi(k+1)}{2(M_{x}+1)}
	\,,
	\quad
	k=0,\ldots,M_{x}-1
	\,,
\end{equation}
where $M_{\tau},M_{x}>0$ define the truncation. The associated discrete inner products are given by
\begin{equation}
	\label{eq:24.05.13_04}
	\langle f, g\rangle_{\tau} = \sum_{j=0}^{M_{\tau}-1}f(\tau_{j})g(\tau_{j})w_{j}
	\,,
	\quad
	w_{j} = \frac{2\pi}{2M_{\tau}+1}
	\,,
\end{equation}
and
\begin{equation}
	\label{eq:24.05.13_05}
	\langle f, g\rangle_{x} = \sum_{k=0}^{M_{x}-1}f(x_{k})g(x_{k})\varpi_{k}
	\,,
	\quad
	\varpi_{k} = \frac{2\pi}{2(M_{x}+1)}
	\,,
\end{equation}
which serve as approximations of their continuous counterparts.

Substituting \eqref{eq:24.05.13_02} into Eq.~\eqref{eq:conformal2} and enforcing the orthogonality of the residual against $\cos(2m+1)\tau\,\sin2(n+1)x$, for $m=0,\ldots,M_{\tau}-1$ and $n=0,\ldots,M_{x}-1$, we get
\begin{multline}
	\label{eq:24.05.13_06}
	\left(\frac{\pi}{2}\right)^{2}\left(-\Omega^{2}(2m+1)^{2}+(2n+1)^{2}\right)\hat{u}_{mn} 
    \\
    + 
	\sum_{j=0}^{\tilde{M}_{\tau}-1}\sum_{k=0}^{\tilde{M}_{x}-1}\frac{\left(u_{M_{\tau},M_{x}}(\tilde{\tau}_{j},\tilde{x}_{k})\right)^{3}}{\sin^{2}{x}}\cos(2m+1)\tilde{\tau}_{j}\,\sin(2n+1)\tilde{x}_{k}\,\tilde{w}_{j}\tilde{\varpi}_{k} = 0
	\,,
\end{multline}
where the last term represents the Fourier coefficients of the nonlinear term, computed in the physical space.

To remove the aliasing errors in \eqref{eq:24.05.13_06} we use the discrete inner products with increased resolution of $\tilde{M}_{\tau}=3M_{\tau}-1$ and $\tilde{M}_{x}=3M_{x}-1$ quadrature rules $(\tilde{\tau}_{j},\tilde{w}_{j})$ and $(\tilde{x}_{k},\tilde{\varpi}_{k})$, ensuring exact integration and preserving equivalence with the Galerkin method. This pseudo-spectral strategy substantially reduces computational cost compared to the traditional Galerkin approach. To improve performance even further, the Jacobian of the equations is evaluated numerically using a finite difference approximation.

The energy corresponding to solution \eqref{eq:24.05.13_02} can be determined by substituting the truncated expansion  into \eqref{eq:25.03.24_05} and evaluating the integral at $\tau = \pi/2$, yielding
\begin{equation}
	\label{eq:24.05.13_08}
	E = \frac{\pi}{4}\Omega^{2} \sum_{n=0}^{M_{x}-1}\left(\sum_{m=0}^{M_{\tau}-1}(-1)^{m}(2m+1)\hat{u}_{mn}\right)^{2}
	\,.
\end{equation}

The problem of finding time-periodic solutions is thereby reduced to solving the nonlinear algebraic system
\begin{equation}
	\label{eq:24.05.12_01}
	\mathbf{F}(\hat{\textbf{u}},\Omega)=0\,,
	\quad
	\mathbf{F}:\mathbb{R}^{M_{\tau}M_{x}}\times \mathbb{R} \rightarrow \mathbb{R}\,,
\end{equation}
where $\hat{\textbf{u}}$ represents the vector of Fourier coefficients in \eqref{eq:24.05.13_02}, and $\mathbf{F}$ denotes the equations stated in \eqref{eq:24.05.13_06}.
To explore the bifurcations along the path of solutions emerging from $(\hat{\textbf{u}},\Omega)=(\textbf{0},N)$, we apply the pseudo-arclength continuation method \cite{Keller.1987}. This strategy allows for efficient solution tracking even when the Jacobian of $\mathbf{F}$ becomes singular. The method is iterated along the solution path, ensuring robustness even for large-amplitude solutions.

\subsection{Evaluation of projection integrals}
\label{sec:EvaluationOfProjectionIntegrals}

The same numerical tools can be used for a quick and reliable computation of the coefficients of the perturbative expansion described in Sec.~\ref{sec:PerturbativeExpansion}.
In particular, we compute the projections appearing in \eqref{eq:25.04.07_01a}, \eqref{eq:25.04.07_01b}, and \eqref{eq:25.04.07_01c} using the Galerkin approach outlined above. To avoid aliasing errors, the employed grids are set large enough to accommodate all modes appearing in the expressions inside the square brackets in \eqref{eq:25.04.07_01} at a given order.
To simplify the implementation, we work with a fixed grid whose size is adjusted to the chosen mode number $N$ and the maximal perturbative order $n_{\textrm{max}}$.
To minimize rounding errors, all calculations are performed with extended precision, typically set to $2n_{\textrm{max}}$ digits.
To verify the correctness of the solution at each order, we monitor the residual and ensure that it remains small and consistent with prescribed precision.
This combined technique allowed us to compute solutions up to $n_{\textrm{max}}=248$ for $N=2$ using the available computing resources. The only restriction was operating system memory, which was quickly exhausted during data output, a current limitation of \textit{Wolfram Mathematica} \cite{Mathematica}. Otherwise, the implementation of the algorithm could have been pushed even further.

\section{Structure of solutions}
\label{sec:StructureOfSolutions}

In this section we discuss the results for $N=2$, i.e., we study the solution family bifurcating from the mode $\cos\tau\,\sin2x$ at the frequency $\Omega=2$. Results for other $N\geq 2$ can be obtained with the same methods and are qualitatively similar.

\subsection{Numerical results}
\label{sec:NumericalResults}

\begin{figure}[t]
    \centering
    \includegraphics[width=0.95\linewidth]{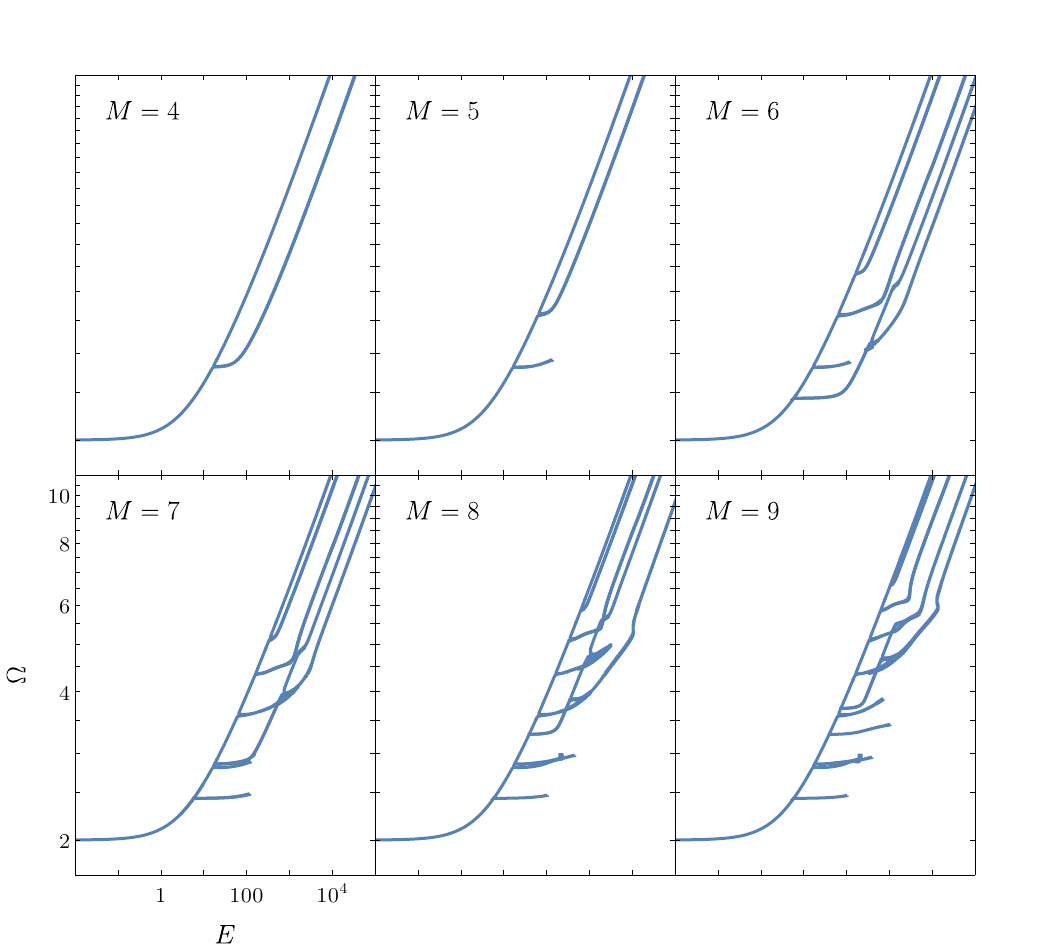}
    \caption{Development of the bifurcation diagram of time-periodic solutions with increasing truncation order $M_{\tau}=M_{x}=M$ in the Galerkin system \eqref{eq:24.05.13_06}. For smaller truncations $M<4$, only the primary trunk bifurcating from the fundamental mode $\cos{\tau}\sin{2x}$ is visible. As the truncation order increases, the diagram becomes progressively more intricate, revealing additional branches and a richer structure. Notice that some of the branches, which initially appear unbounded, become bounded as $M$ increases.}
    \label{fig:Galerkin}
\end{figure}

\begin{figure}[t]
    \centering
    \includegraphics[width=0.95\linewidth]{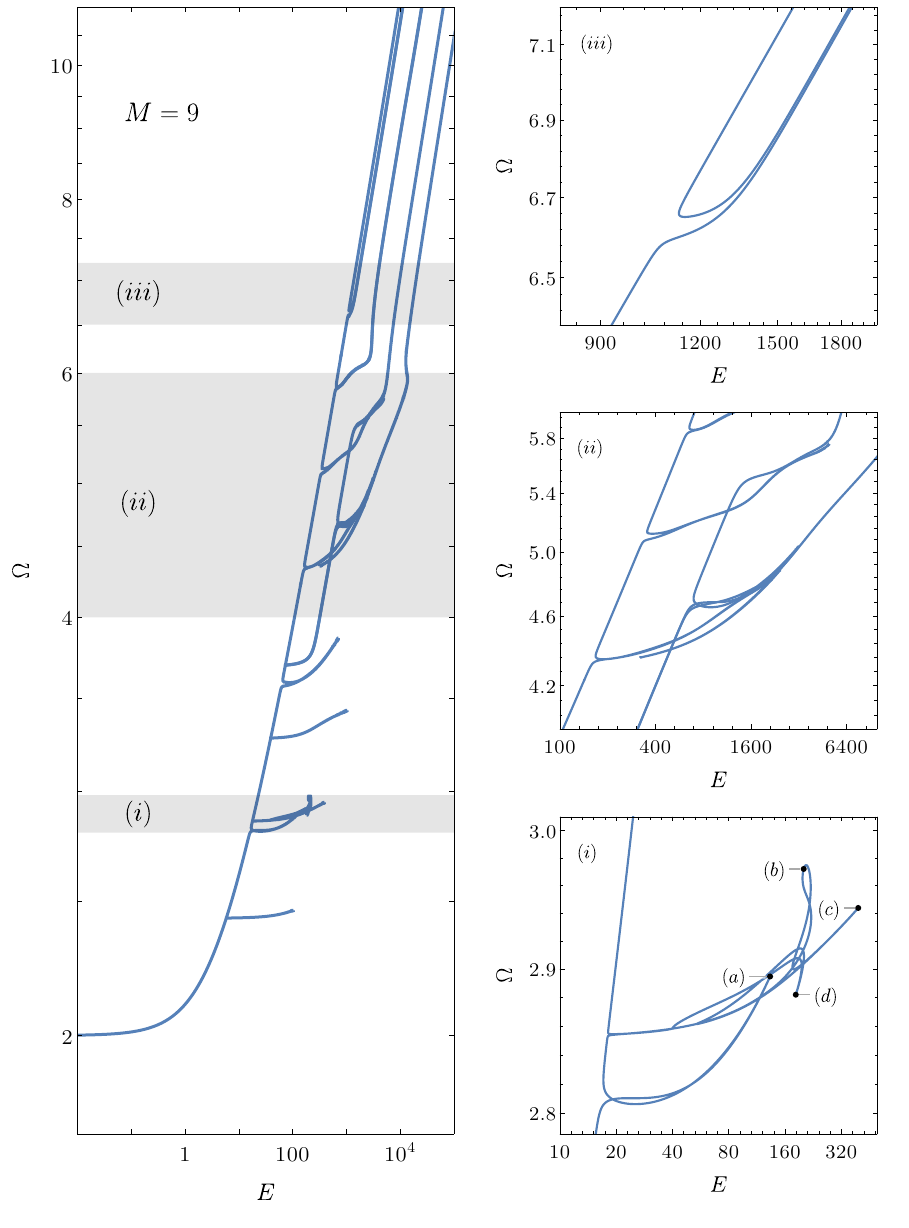}
    \caption{Details of the $\Omega-E$ diagram for mode truncation $M_{\tau}=M_{x}=M=9$, along with several zoom-in panels. See Fig.~\ref{fig:GalerkinModes} for the behaviour of the mode amplitudes associated with the branches shown in panel $(i)$.}
    \label{fig:GalerkinZoom}
\end{figure}

The Galerkin-based approach outlined in the previous section reveals similar solution structures to the ones discovered for the cubic wave equation in \cite{Ficek.2024A}. Specifically, we observe the emergence of a trunk extending indefinitely, accompanied by a network of branches located at discrete points along the trunk. However, in the case at hand, we observe that some branches also extend indefinitely, continuing to arbitrarily high energies and frequencies. Nevertheless, as the truncation order of the Galerkin system increases, these branches become bounded, see Fig.~\ref{fig:Galerkin} and \ref{fig:GalerkinZoom}. With further increases in truncation, additional structures may emerge from the branches, but their overall shapes remain stable. In particular, the ends of the branches remain connected to the rescaled trunks bifurcating from the corresponding modes (see discussion below).

Compared to the cubic wave equation, where the bifurcation structure is relatively simple, likely attributed to the scaling symmetry, the conformal case exhibits significantly more intricate patterns in the energy-frequency diagram as illustrated in Fig.~\ref{fig:GalerkinZoom}. The solution paths twist and fold in a complex manner, reflecting the richer nonlinear interactions present in the conformally invariant setting and the lack of scaling symmetry.

\begin{figure}[t]
    \centering
    \includegraphics[width=0.55\linewidth]{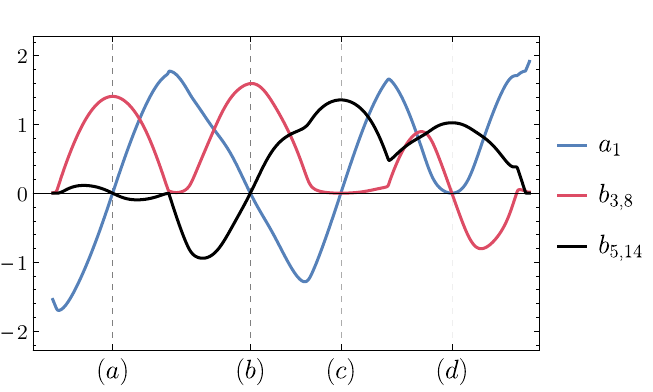}
    \caption{The amplitudes of modes $\cos\tau\,\sin 2x$ (blue), $\cos3\tau\,\sin 8x$ (red), $\cos5\tau\,\sin 14x$ (black) along the part of the diagram shown in Fig.~\ref{fig:GalerkinZoom}. The marked locations $(a),\ldots, (d)$ correspond to the points indicated on the $\Omega-E$ diagram. We observe the fundamental mode to decrease in amplitude when it enters the branches and vanish at the bifurcation points $(a),\ldots, (d)$, that also lie on the appropriate rescalings of trunks bifurcating from zero at other frequencies.}
    \label{fig:GalerkinModes}
\end{figure}

\begin{figure}
    \centering
    \includegraphics[width=0.48\linewidth]{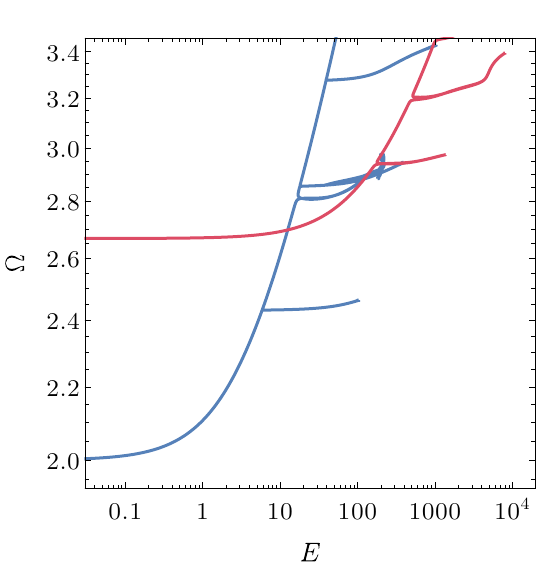}
    \hspace{2ex}
    \includegraphics[width=0.48\linewidth]{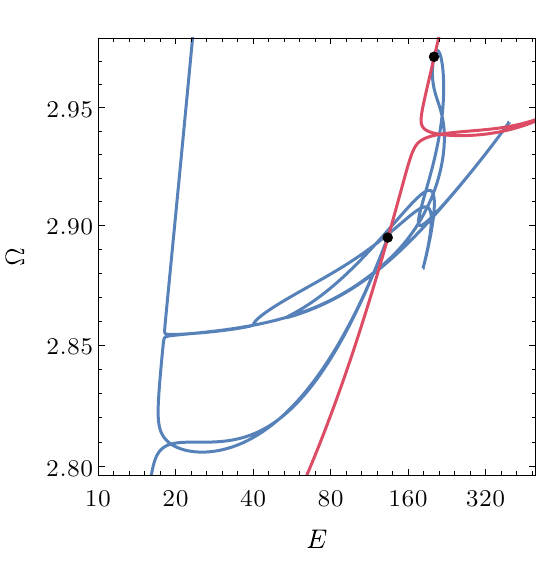}
    \caption{
    The bifurcation structure near the locations of the branches dominated by the modes $\cos{3\tau}\,\sin{8x}$ and $\cos{5\tau}\,\sin{14x}$. The solution bifurcating from $\cos{\tau}\,\sin{2x}$, as shown in Figs.~\ref{fig:Galerkin} and \ref{fig:GalerkinZoom}, is plotted in blue. The red curve corresponds to the solution bifurcating from $\cos{\tau}\,\sin{8x}$ which has been rescaled as described in the text and intersects the bifurcation points indicated by black dots (points $(a)$ and $(b)$ in Fig.~\ref{fig:GalerkinZoom}). The solution bifurcating from $\cos{5\tau}\,\sin{14x}$, when appropriately rescaled, would match the remaining bifurcation points $(c)$ and $(d)$ (for clarity, this curve is not shown).
    }
    \label{fig:Branch_38}
\end{figure}

As in the case of the cubic wave equation, the observed structures exhibit a characteristic interplay between fundamental and higher harmonics: upon entering a bounded branch, higher modes become dominant, while the fundamental mode is suppressed, only to reverse as the branch returns to the trunk, see panel $(i)$ in Fig.~\ref{fig:GalerkinZoom} and the evolution of mode amplitudes along this part of the plot shown in Fig.~\ref{fig:GalerkinModes}. The points $(a),\ldots, (d)$ marked there are bifurcation points connecting the investigated solution family (bifurcating from the mode $\cos\tau\,\sin 2x$) with appropriate rescalings of other solution families. For example, points $(a)$ and $(b)$ lie on a trunk of solutions bifurcating from the mode $\cos3\tau\,\sin8x$ at $\Omega=8/3$, see Fig.~\ref{fig:Branch_38}. This trunk can be obtained from the one bifurcating from $\cos\tau\,\sin8x$ at $\Omega=8$ (hence, belonging to the family of solutions investigated in this work) by a simple rescaling $\tau\to3\tau$, $\Omega\to\Omega/3$. 
This intricate branching pattern, governed by nonlinear mode interactions, suggests a universal mechanism underlying the formation of a complex structure of time-periodic solutions. 

We remark that convergence of the Galerkin approximation was verified by monitoring the $L^{2}$ norm of the residual as a function of the truncation size $M = M_{\tau} = M_{x}$. In all cases considered, we observe exponential (spectral) convergence toward machine precision, with minor deviations associated with the emergence of branches and the spectral structure of the solutions, cf. \cite{Ficek.2024A} for details.

\subsection{Reducible systems}
\label{sec:ReducibleSystems}

We can employ the method based on the reducible systems discussed in \cite{Ficek.2024B} as an alternative, independent approach to the study of the locations of the branches. Reducible systems are collections of modes that are coupled with each other in a minimal way\footnote{See discussion in Sec.~3 of \cite{Ficek.2024B}.} leading to exceptionally simple Galerkin systems of equations. 

The overall shape of the trunk can be approximated by a one mode system spanned by $A\cos\tau\,\sin Nx$. Then the Galerkin scheme leads to a single equation
\begin{align}
A \left(3N A^2-4\Omega^2+4N^2 \right)=0\,
\end{align}
with a non-trivial solution given by
\begin{align}
    A=2\sqrt{\frac{\Omega^2-N^2}{3N}}\,.
\end{align}
It exists for $\Omega>N$ and has energy $E=\pi\,\Omega^2 A^2/4$. 

The observed structures emerging from the trunk can be recreated using reducible systems spanned by two modes: $A\cos\tau\,\sin Nx$ and $B\cos m\tau \, \sin n x$, where $m$ is odd and $n$ has the same parity as $N$. For these modes to compose a reducible system, it must hold either $m\geq 5$ and $n\geq N$, or $m=3$ and $n\geq 3N$. Then the Galerkin system of equations spanned by them is
\begin{equation}
    \label{eq:2modes_reducible}
    \left\{
    \begin{aligned}
        A \left(3N A^2+6N B^2-4\Omega^2+4N^2 \right)&=0\,,
        \\
        B \left(6N A^2+ 3n B^2-4\Omega^2 m^2+4n^2\right)&=0\,.
    \end{aligned}
    \right.
\end{equation}
The energy of solutions to such systems is given by
\begin{align}
    E=\frac{\pi}{4}\Omega^2\left(A^2+  m^2 B^2\right)\,.
\end{align}
For $n\neq 4N$ solutions to \eqref{eq:2modes_reducible} are formally given by 
\begin{align}
    A=2\sqrt{\frac{(2n-N)nN+(n-2m^2 N)\Omega^2}{3(n-4N)N}}\,,\qquad
    B=2\sqrt{\frac{2N^2-n^2+(m^2-2)\Omega^2}{3(n-4N)}}\,.
\end{align}
The analysis of the expressions under the roots leads us to the following conditions under which $A$ and $B$ are real. For $m=3$, $n$ must be inside the interval $3N+2\leq n\leq 4N$. Since $n$ has the same parity as $N$, $n$ is every second integer inside this interval, and the upper boundary is included only for even $N$. For $m>3$ the condition is simply $n=mN+2k$ where $k\in\mathbb{N}$. For fixed $N$ any pair of $(m,n)$ satisfying these conditions gives a two-modes solution branch bifurcating from the trunk at $\Omega=\sqrt{(n^2-2N^2)/(m^2-2)}$. Similarly as in \cite{Ficek.2024B}, it lets us suspect that branches densely populate the trunk, which stays in agreement with the previous results \cite{Berti.2024,Silimbani.2024, Chatzikaleas.2020oh, Chatzikaleas.2024}. However, these branches represent the new class of solutions that has not been considered earlier.

When $n$ and $N$ are both even, we need to additionally consider the case $n = 4N$. It leads to an overdetermined linear system that has solutions only for 
$\Omega=N\sqrt{14/(m^2-2)}$. Then the amplitudes of two-mode solutions must satisfy
\begin{align}
    3A^2+6B^2+4N-\frac{56N}{m^2-2}=0\,.
\end{align}
Since it may hold only for $m=3$, we get a single additional branch bifurcating from the trunk at the frequency $\Omega=\sqrt{2}N$.

As an example, for $N=2$, in Tab.~\ref{tab:reducible} we present pairs of numbers $(m,n)$ such that reducible systems spanned by $A\cos\tau\,\sin 2x$ and $B\cos m\tau \, \sin n x$ give branches bifurcating from the main trunk. By selecting from this table all entries with modes $(m,n)$ such that $m\leq 2N-1$ and $n\leq 2N$ and calculating frequencies of the aforementioned bifurcation points, we get numerical values that agree with the locations of branches for respective truncation $M$ on Fig.~\ref{fig:Galerkin}. The agreement is better for the lower branches, i.e., branches with the frequency $\Omega$ close to $\Omega=2$.
\begin{table}
    \begin{center}
    \begin{tabular}{ccccccccc}
    \colorbox{verylightgray}{$(3,8)$} & $(3,10)$ & $(3,12)$ & $(3,14)$ & $(3,16)$ & $(3,18)$ & $(3,20)$ & $(3,22)$ & $\hdots$ \\
     & & \colorbox{verylightgray}{$(5,12)$} & $(5,14)$ & $(5,16)$ & $(5,18)$ & $(5,20)$ & $(5,22)$ & $\hdots$ \\
      & & & & \colorbox{verylightgray}{$(7,16)$} & $(7,18)$ & $(7,20)$ & $(7,22)$ & $\hdots$ \\
        & & & & & & \colorbox{verylightgray}{$(9,20)$} & $(9,22)$ & $\hdots$\\
            & & & & & & & & $\ddots$
    \end{tabular}
    \end{center}
    \caption{Mode numbers $(m,n)$ which lead to branches spanned by $A\cos\tau\,\sin 2x$ and $B\cos m\tau \, \sin n x$ for solutions bifurcating from mode $N=2$. In Sec.~\ref{sec:PadeApproximation} we expand the coefficients of highlighted modes into Pad\'{e} approximants in $\varepsilon$, then the values of $\Omega$ for which these approximants have poles are plotted in Fig.~\ref{fig:spectrum}.
    }
    \label{tab:reducible}
\end{table}

Let us point out that in this work the reducible systems approach is employed only as a tool for establishing approximate locations of the branches emerging from the trunk. Although it succeeds in this regard, it does not recreate the higher-energetic parts of the branches very well, in contrast to what was observed in \cite{Ficek.2024B}. This difference likely comes from  the additional interplays between the modes, introduced by the much richer structure of interaction coefficients (see Eq.\ \eqref{eq:decomposition_main} and discussion there) that cannot be approximated by such simple systems.

\subsection{Pad\'{e} approximation}
\label{sec:PadeApproximation}
\begin{figure}
    \centering
    \includegraphics[width=0.48\linewidth]{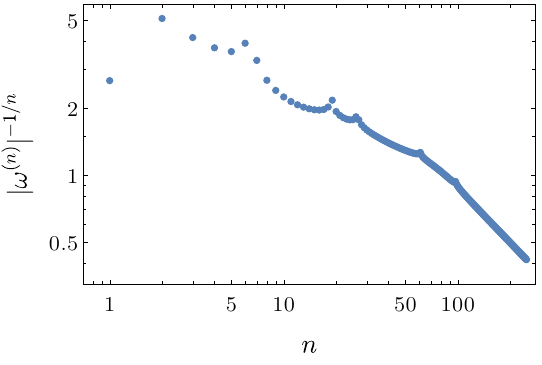}
    \hspace{2ex}
    \includegraphics[width=0.48\linewidth]{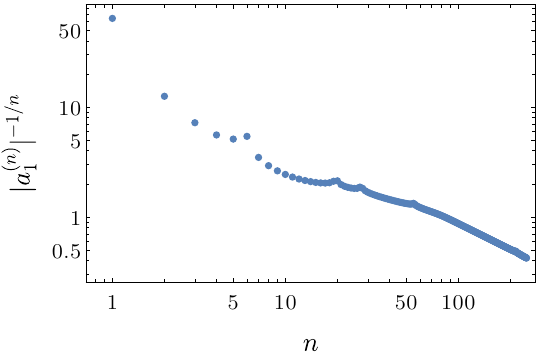}
    \caption{Behaviour of coefficients in perturbative expansions of $\Omega$ and $a_1$ for solution bifurcating from $\cos\tau\,\sin 2x$. It suggests that the radius of convergence for both series is zero.}
    \label{fig:convergence_radius}
\end{figure}

Let us write the solution $u$ to Eq.~\eqref{eq:conformal2} as a Fourier series
\begin{align}
    u(\tau,x)=\sum_{m=0}^\infty a_{2m+1} \cos(2m+1) \tau\,\sin (2m+1)Nx +\sum_{j=0}^\infty \sideset{}{'}\sum_{k=1}^{\infty}b_{2j+1,k}\cos(2j+1)\tau\,\sin kx\,.
\end{align}
We can treat every coefficient $a_{2m+1}$ and $b_{2j+1,k}$ as a formal power series in $\varepsilon$ with terms that can be constructed up to arbitrary order via the algorithm discussed in Sec.~\ref{sec:PerturbativeExpansion}. The same observation applies to the frequency $\Omega$. We suspect that such series have zero radius of convergence, see Fig.~\ref{fig:convergence_radius} and discussion in \cite{arnold1988dynamical}. Nevertheless, as we show in this section, the Pad\'{e} approximations that come from them encode a lot of information regarding the structure of time-periodic solutions.

\begin{figure}[t]
    \centering
    \includegraphics[width=0.99\linewidth]{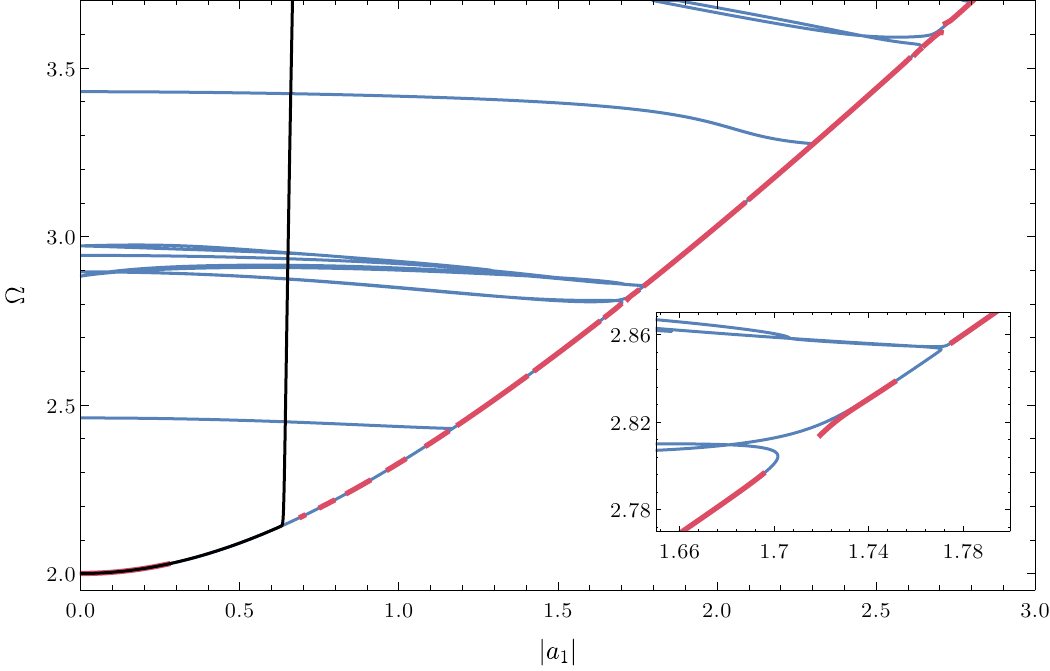}
    \caption{Curves in $(|a_1|,\Omega)$ space formed by the solutions to the numerical scheme with $M_\tau=M_x=M=9$ (blue), the power series generated by the perturbative expansion up to the 248-th order (black), and the corresponding Pad\'{e} approximation of type $[124/124]$ (red). Note that the agreement between the numerical (blue) and Pad\'{e} (red) curves extends well beyond the perturbative regime along the trunk, except in the vicinity of the branches; for details see the inset.
    }
    \label{fig:pade1}
\end{figure}

For any fixed positive integer $n$, we can consider a power series for a Fourier coefficient, e.g, $a_1$, truncated at order $2n$. Since it is a polynomial in $\varepsilon$ of order $2n$, we can use it to construct a Pad\'{e} approximant of type $[n/n]$. The same procedure can be done for a truncation of a power series for $\Omega$. Thus, we get two expressions parametrised by a common parameter $\varepsilon$. We present the resulting curve $(|a_1|,\Omega)$ for the approximants of type $[124/124]$ in red in Fig.~\ref{fig:pade1}. From this curve, we have removed neighbourhoods of the poles of the Pad\'{e} approximants. One can see a very strong improvement compared to the black curve coming from the original series for $a_1$ and $\Omega$, which breaks down near the value of $\varepsilon$ suggested by Fig.~\ref{fig:convergence_radius}. Similar analysis can be made for any other coefficient $a_{2m+1}$ or $b_{2j+1,k}$, leading to similar observations.

\begin{figure}[t]
    \centering
    \includegraphics[width=0.99\linewidth]{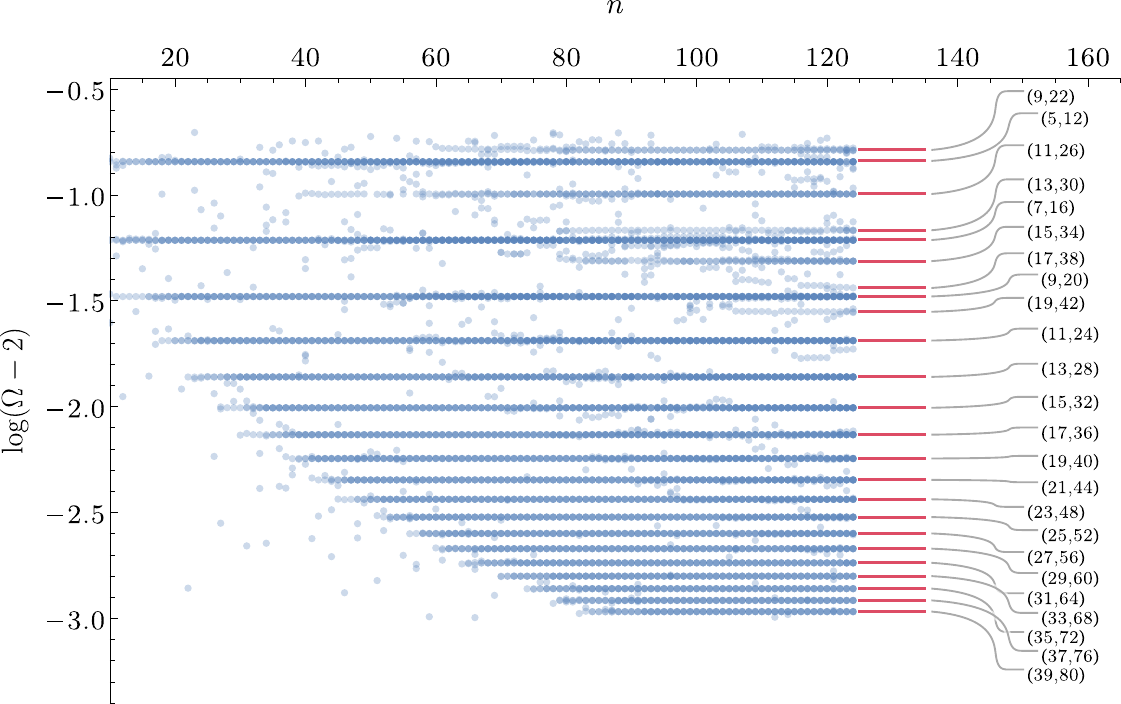}
    \caption{Frequency $\Omega$ evaluated at the poles of $[n/n]$ Pad\'{e} approximants for $b_{3,8}$, $b_{5,12}$, $b_{7,16}$, and $b_{9,20}$ (see Tab.~\ref{tab:reducible}) of the solution bifurcating from $\cos\tau\,\sin 2x$. The red lines indicate selected frequencies at which two-modes solutions, spanned by $\cos\tau\,\sin 2x$ and $\cos m\tau\,\sin nx$, with $(m,n)$ denoted by the associated indices, bifurcate from the trunk, according to the reducible systems approach.}
    \label{fig:spectrum}
\end{figure}

We have mentioned that the red curve in Fig.~\ref{fig:pade1} does not contain neighbourhoods of the poles of the Pad\'{e} approximants. One can observe that holes resulting from this omission are often in the vicinity of the branches. This observation can be explained by the fact that the beginnings of the branches constitute places along the solution curve where the amplitudes undergo strong changes. This type of rapid change can be mimicked by the Pad\'{e} approximant having a pole in its vicinity. It leads us to the following idea of predicting the positions of the branches using the poles of the Pad\'{e} approximants. 
For any fixed coefficient $a_{2m+1}$ or $b_{2j+1,k}$ we can construct Pad\'{e} approximants of the type $[n/n]$.
(We note that no significant improvement was observed when using non-diagonal Pad\'{e} approximants; therefore, for concreteness, we present the results for approximants of the type $[n/n]$.)
Next, we find poles of such approximants by simply looking for real roots of the polynomial of order $n$ in its denominator. Then, these roots can be plugged into a Pad\'{e} approximant of $\Omega$ to give us frequencies at which one may expect a branch. As an example, in Fig.~\ref{fig:spectrum} we show the results of such a procedure for coefficients $b_{3,8}$, $b_{5,12}$, $b_{7,16}$, and $b_{9,20}$ of the solution bifurcating from $\cos\tau\,\sin 2x$. One can see that in addition to the spurious poles that are natural in this type of technique \cite{Baker.1996}, some clear patterns are easily visible. Locations of these patterns are then compared with the frequencies of branches predicted via the reducible systems approach, giving a striking agreement. We also observe that branches associated with modes in the form $(m, 2m+2)$ with $m$ being an odd number (the `leftmost' modes in the extension of Tab.~\ref{tab:reducible}) appear for much lower orders of Pad\'{e} approximants than other modes. Let us note that an analogous analysis done for different choices of coefficients would result in a similar plot, showing various horizontal patterns. Such a plot for all coefficients would reproduce more and more branches as the maximal order $n$ is increased, recreating in the limit $n\to\infty$ their whole infinite structure.

Furthermore, we have confirmed that the structure shown in Fig.~\ref{fig:spectrum} is robust under perturbations of the expansion coefficients. Upon introducing substantial noise to the output of the perturbative scheme--significantly larger than the numerical precision used--we observed the spurious shift ever so slightly, while the overall structure remained intact and consistent with the predictions of the reducible systems. This indicates that the spurious poles are not artifacts of finite numerical precision but intrinsic features of the Pad\'{e} approximation itself.
It is plausible that these scattered poles are precursors to branches that become visible only at higher orders in the perturbative expansion.

\appendix

\vspace{8ex}

\section{Interaction coefficients}
\label{sec:InteractionCoefficients}
\begin{figure}
    \centering
    \includegraphics[width=0.7\linewidth]{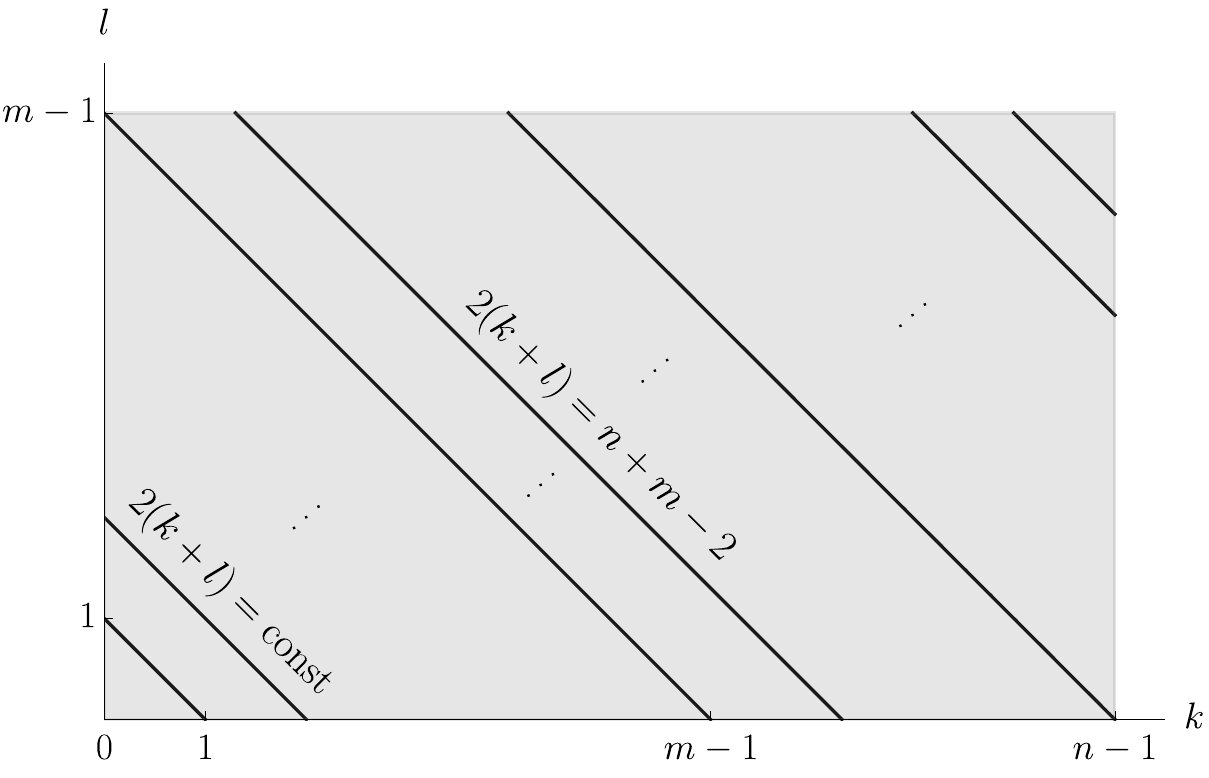}
    \caption{A diagram showing the structure of elements inside the sum \eqref{eq:polynomialP_supp}.}
    \label{fig:polynomialP_supp}
\end{figure}

In this section we show that for $j,k,l\in\mathbb{N}$ it holds
\begin{align}
    \label{eq:decomposition}
    \frac{\sin jx\, \sin k x\, \sin l x }{\sin^2 x}=\sum_{m=1}^{j+k+l-2} S_{jklm}\, \sin mx,
\end{align}
where interaction coefficients $S_{jklm}$ can be written with the use of functions
\begin{align}
    \label{eq:min}
    \mathfrak{m}(j,k)=\frac{1}{2} \sgn j\, \sgn k\, \min(|j|,|k|)
\end{align}
as 
\begin{equation}
    \label{eq:S}
    S_{jklm} = \begin{cases}
        \mathfrak{m}(j+k-l,m)+\mathfrak{m} (j-k+l,m)+\mathfrak{m} (-j+k+l,m)- \mathfrak{m} (j+k+l,m)\\
          \qquad\qquad\qquad\qquad\qquad\qquad\qquad\qquad\qquad\qquad \mbox{ if } j+k+l+m  \mbox{ even}\,,\\
          \\
        0 \qquad\qquad\qquad\qquad\qquad\qquad\qquad\qquad\qquad\quad\ \,\mbox{ if } j+k+l+m  \mbox{ odd}\,.
    \end{cases}
\end{equation}
Thus, we want to calculate
\begin{align*}
    S_{jklm}=\frac{2}{\pi}\int_0^\pi \frac{\sin jx\, \sin k x\, \sin l x\, \sin m x }{\sin^2 x}\,\dd x\,.
\end{align*}
The simple formula
\begin{multline*}
    \sin j x \, \sin k x\, \sin l x
    \\
    = \frac{1}{4}\left[\sin(j+k-l)x+\sin(j-k+l)x+\sin(-j+k+l)x-\sin(j+k+l)x\right]
\end{multline*}
lets us rewrite this expression as
\begin{align}
    \label{eq:SCCCC_supp}
    S_{jklm}= \frac{1}{4}\left(C_{j+k-l,m}+C_{j-k+l,m}+C_{-j+k+l,m}-C_{j+k+l,m}\right)\,,
\end{align}
where 
\begin{align}
    C_{nm}=\frac{2}{\pi}\int_0^\pi \frac{\sin n x\, \sin m x}{\sin^2 x}\,\dd x\,,\qquad m,n\in\mathbb{Z}\,.
\end{align}
To find $C_{nm}$, let us initially assume that $n,m>0$. We introduce a complex variable $t=e^{ix}$. Then, we have the following equalities
\begin{multline*}
    C_{nm}=\frac{1}{\pi}\int_0^{2\pi} \frac{\sin n x\, \sin m x}{\sin^2 x}\,\dd x=\frac{1}{\pi}\int_0^{2\pi}\frac{\left(e^{inx}-e^{-inx}\right)\left(e^{imx}-e^{-imx}\right)}{\left(e^{ix}-e^{-ix}\right)^2}\,\dd x
    \\
    = \frac{1}{i\pi}\oint_\gamma \frac{\left(t^n-t^{-n}\right)\left(t^m-t^{-m}\right)}{\left(t-t^{-1}\right)^2 t}\,\dd t
    = \frac{1}{i\pi}\oint_\gamma t^{-n-m+1} \frac{\left(1-t^{2n}\right)\left(1-t^{2m}\right)}{\left(1-t^{2}\right)^2}\,\dd t
    \\
    = \frac{1}{i\pi}\oint_\gamma t^{-n-m+1} \frac{\left(1+t+\ldots+t^{2n-1}\right)\left(1+t+\ldots+t^{2m-1}\right)}{\left(1+t\right)^2}\,\dd t\,,
\end{multline*}
where $\gamma$ is the unit circle in the complex plane. One can gather terms inside the brackets to get
\begin{align*}
    1+t+\ldots+t^{2n-1} = (1+t)\left(1+t^2+\ldots+t^{2n-2}\right)\,.
\end{align*}
As a result, we obtain
\begin{align*}
    C_{nm}=& 
    \frac{1}{\pi i}\oint_\gamma t^{-n-m+1} P(t)\,\dd t= 2 \,\mbox{Res}\left(t^{-n-m+1}P(t),0\right)\,.
\end{align*}
where we have defined the polynomial $P$ as
\begin{align}
    \label{eq:polynomialP_supp}
    P(t)=\left(1+t^2+\ldots+t^{2n-2}\right)\left(1+t^2+\ldots+t^{2m-2}\right)=\sum_{k=0}^{n-1} \sum_{l=0}^{m-1}t^{2(k+l)}\,.
\end{align}
The residue is given just by the coefficient next to $t^{n+m-2}$ term in $P$. Obviously, if $n$ and $m$ have different parities, there is no such term, so the result is zero. Therefore, we focus on the case when they both are either even or odd. Without the loss of generality, let us assume that $m\leq n$. Then, the graphical representation of the sum in \eqref{eq:polynomialP_supp}, see Fig.~\ref{fig:polynomialP_supp}, tells us that terms of orders between $2m-2$ and $2n-2$ have multiplicity $m$. Since $n+m-2$ belongs to this range, the residuum is equal to $m$. It means that for $n,m>0$, $C_{nm}=2\min(n,m)$. A simple case-by-case analysis lets us conclude that when $n$ and $m$ are of the same parity, then
\begin{align*}
    C_{nm}=2\,\sgn \, m\; \sgn\, n \: \min(|n|,|m|)\, ,
\end{align*}
while $C_{nm}=0$ in the opposite case. Plugging this result into Eq.~\eqref{eq:SCCCC_supp} gives us Eq.~\eqref{eq:S}.

In the end, let us point out that from Eq.~\eqref{eq:S} one quickly gets that $S_{jklm}=0$ when $m\geq j+k+l-1$. Thus, it is sufficient to set the upper limit of summation in \eqref{eq:decomposition} to $j+k+l-2$. In addition, the symmetry of $S_{jklm}$ in all indices leads to the conclusion that $S_{jklm}=0$ if the sum of any three indices is equal to or larger than the fourth one. We also note two particularly useful formulas:
\begin{align}
    S_{nnnn}=n\,, \qquad S_{nnkk}=\min(n,k)\,.
\end{align}

\printbibliography

\end{document}